# Determining Maximal Reference Set in Data Envelopment Analysis

I. Roshdi $^{a,*}$, I. Van de Woestyne $^{b}$, M. Davtalab–Olyaie $^{c}$

$^{a}$*Department of Mathematics, Science and Research Branch, Islamic Azad University, Tehran, Iran*
$^{b}$*Hogeschool − Universiteit Brussel, Warmoesberg 26, B − 1000 Brussel, Belgium*
$^{c}$*Faculty of Mathematical Sciences and Computer, Kharazmi University, Tehran, Iran*

*Abstract*

In data envelopment analysis (DEA), the occurrence of multiple reference sets is a crucial issue in identifying all the reference DMUs to a given decision making unit (DMU). To resolve this difficulty, we introduce the useful notion of *maximal reference set* (MRS) which contains all the reference DMUs. Based on both primal and dual formulations, we then propose a new mixed integer linear programming (MILP) based approach for determining the MRS. The proposed approach is more general than the existing ones, and has several desirable properties: (i) having a unified formulation, (ii) adaptable with different DEA models, and (iii) compatible with both primal and dual forms of the DEA models. Furthermore, our approach is made computationally effective by establishing a LP-based procedure to treat with the primal-based MILP.

***Keywords:*** *Data envelopment analysis, Efficiency, Reference set, Mixed integer linear programming.*

## 1. Introduction

DEA, introduced by Charnes et al. [6] and subsequently extended by Banker et al. [1], is a linear programming (LP) based methodology for measuring the relative efficiency of a set of homogeneous DMUs with multiple inputs and multiple outputs. In addition to the efficiency score, DEA provides a set of efficient projection points for the DMU being evaluated. The coordinates of each projection can be interpreted as the "target" levels of operation of inputs and outputs that indicate how the assessed DMU can be improved to perform efficiently. A projection point (if it is not an observed DMU) is represented as a convex combination of the efficient DMUs in a set known as a *reference set* (RS). Those efficient

* Corresponding author, Email: i.roshdi@gmail.com

DMUs forming the RS are known as the *reference DMUs* to which the assessed DMU is directly compared to obtain its efficiency. Identification of the RS plays a crucial role in many areas of the DEA including:

- sensitivity and stability analysis (Hibiki and Sueyoshi [10], Boljuncic [3]),
- identifying the status of returns to scale (Sueyoshi and Sekitani [19], Krivonozhko et al. [13,14]),
- ranking efficient DMUs (Jahanshahloo et al. [11]),
- benchmarking and target setting (Bergendahl [2], Camanho and Dyson [4]).

As such, determination of the RS has recently received considerable attention in DEA literature. In a seminal work and based on the strong complementary slackness conditions between the primal and dual of the BCC model, Sueyoshi and Sekitani [19] provided some characterization of the RS. Furthermore, they developed an integrated primal-dual formulation for the RS determination. However, as demonstrated by Krivonozhko et al. [15], "not only in their approach the computational burden increases significantly, but also it seems that the basic matrices may be inherently ill-conditioned, leading to erroneous results." Unlike the single-stage approach in [19], Jahanshahloo et al. [12] proposed a computational multi-stage procedure to identify the RS based on the multiplier (dual) form of the BCC model and tried to reduce the size of the problem in each stage. However, their approach is not applicable for a DEA model with nonlinear objective such as the Russell and enhanced Russell measures [17, 18].

Corresponding to the BCC model, we first introduce the useful notion of MRS containing all the reference DMUs to the evaluated DMU. Then, pertaining to both primal and dual forms of the BCC model, we develop an MILP-based approach to determine MRS. Since this approach is compatible with both primal and dual formulation of DEA models, it can be implemented through primal or dual whichever is desired. Moreover, as we will extend in Subsection 3.3, our approach can be easily applied to all kind of DEA models, and thus is more general than the existing ones. The proposed approach involves solving a unique MILP program thereby it presents a unified formulation to the identification of MRS. To enhance the computational capability of the primal-based MILP program, integrating the concepts of LP and DEA, we establish an LP-based procedure. Moreover, we reveal a linkage between our dual-based MILP and the procedure proposed by Jahanshahloo et al. [12]; indeed, the later can be interpreted as a computational treatment for the former.

## 2. Preliminaries

We assume that we are dealing with $n$ observed DMUs; each uses $m$ inputs to produce $s$ outputs. We denote by $x_j = (x_{1j}, \dots, x_{mj})^T \in \text{IR}_{\geq 0}^m$ and $y_j = (y_{1j}, \dots, y_{sj})^T \in \text{IR}_{\geq 0}^s$ the input and output vectors of $\text{DMU}_j, j \in J = \{1, \dots, n\}$. We assume that $x_j$ and $y_j$ are nonzero vectors for each $j \in J$.

The production possibility set (PPS) $T$ is defined as the set of vectors $(x, y)$ such that the output vector $y$ can be produced by the input vector $x$, i.e.,

$$T = \{(x, y) | x \in \mathrm{IR}_{\geq 0}^m \text{ can produce } y \in \mathrm{IR}_{\geq 0}^s\}. \tag{1}$$

Banker et al. [1] have deduced the following PPS, considering some postulates. This set is denoted by $T_v$, regarding the prevalence of variable returns to scale assumption:

$$T_v = \left\{ (x, y) \in \mathrm{IR}_{\geq 0}^{m+s} \,\middle|\, x \geq \sum_{j \in J} \lambda_j x_j, y \leq \sum_{j \in J} \lambda_j y_j, \sum_{j \in J} \lambda_j = 1, \lambda_j \geq 0, j \in J \right\}. \tag{2}$$

Relative to this PPS, the efficiency of $\mathrm{DMU}_o, o \in J$, can be measured via the following LP, which is called the input-oriented *envelopment* form of the BCC model:

$$\begin{aligned}
\text{Min } & \theta \\
\text{s.t. } & \sum_{j \in J} \lambda_j x_{ij} \leq \theta x_{io}, \quad i = 1, \ldots, m, \\
& \sum_{j \in J} \lambda_j y_{rj} \geq y_{ro}, \quad r = 1, \ldots, s, \\
& \sum_{j \in J} \lambda_j = 1, \\
& \lambda_j \geq 0, j \in J.
\end{aligned} \tag{3}$$

Model (3) is *radial* as it deals only with the proportional decrease in inputs and, thus, fails to account for the additional sources of inefficiency caused by the non-zero slacks. To remove this difficulty, the following LP problem is solved in the second phase [9]:

$$\begin{aligned}
\text{Max } & \sum_{i=1}^{m} s_i^- + \sum_{r=1}^{s} s_r^+ \\
\text{s.t. } & \sum_{j \in J} \lambda_j x_{ij} + s_i^- = \theta^* x_{io}, \quad i = 1, \ldots, m, \\
& \sum_{j \in J} \lambda_j y_{rj} - s_r^+ = y_{ro}, \quad r = 1, \ldots, s, \\
& \sum_{j \in J} \lambda_j = 1, \\
& \lambda_j \geq 0, j \in J, s_i^- \geq 0, s_r^+ \geq 0, \forall i, r.
\end{aligned} \tag{4}$$

where $\theta^*$ is the optimal objective of (3).

**Definition 1.** $\mathrm{DMU}_o$ is said to be BCC-efficient if and only if it is *radial efficient* (i.e., $\theta^* = 1$) and the optimal objective of (4) is equal to zero.

## 3. Characterizing and Finding the Maximal Reference Set

Let us turn to Model (3). We define the reference set of $DMU_o$ as follows:

**Definition 2.** (*Reference Set*) Let $\lambda^*$ be an optimal solution of (3) corresponding to a given $DMU_o$. We define the set of DMUs corresponding to positive $\lambda_j^*$ as the reference set (RS) to $DMU_o$, and denote it by $R_o$ as

$$R_o = \{DMU_j | \lambda_j^* > 0 \}. \tag{5}$$

Each member of the set $R_o$ is called a *reference DMU* to $DMU_o$.

Because of the occurrence of alternative optimal solutions for variable $\lambda$ in (3), the RS is not unique. However, we can go further and determine such an RS containing *all* the RSs to $DMU_o$.

**Definition 3.** (*Maximal Reference Set*) We define the set containing *all* the RSs to $DMU_o$ as the maximal reference set (MRS) for $DMU_o$ and denote it by $R_o^M$.

A desirable property of $R_o^M$ is that it contains all reference DMUs to $DMU_o$. It can be easily verified that each reference DMU is radial efficient [12]. From the target-setting point of view, it is valuable to identify Pareto-efficient (BCC-Efficient) peers or reference DMUs to a given inefficient DMU. Therefore, we filter out the BCC-inefficient reference DMUs, and define $E_o^M$, as the set comprising all the BCC-efficient reference DMUs to $DMU_o$.

To provide a visual representation of the above definitions, we give the following numerical example.

**Example 1.** Consider a single-input- single-output technology comprising seven DMUs whose data are demonstrated in Table 1. The empirical structure of the set $T_v$ is depicted in Figure 1.

**Table1.** Input and output data for Example 1

|  | DMU₁ | DMU₂ | DMU₃ | DMU₄ | DMU₅ | DMU₆ | DMU₇ |
|---|---|---|---|---|---|---|---|
| **Input (I)** | 1.00 | 2.00 | 3.00 | 4.00 | 6.00 | 5.00 | 4.00 |
| **Output (O)** | 1.00 | 3.00 | 4.00 | 5.00 | 6.00 | 4.00 | 3.00 |

In evaluating the inefficient unit $DMU_6$ via Model (3), the corresponding projection point is $DMU_3$ and its possible RSs are $\{DMU_3\}$ and $\{DMU_2, DMU_4\}$. Thus, $R_6^M = \{DMU_2, DMU_3, DMU_4\}$ is the MRS to $DMU_6$. Since all members of the set $R_6^M$ are BCC-efficient, $R_6^M = E_6^M$.

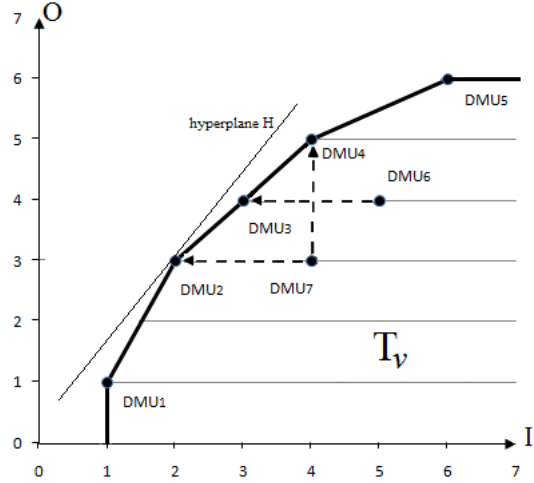

**Figure 1.** The PPS and multiple reference sets

We develop an approach to determine the MRS based upon both primal and dual from of the BCC model in details.

## 3.1. Primal Approach

To determine $R_o^M$, among all the optimal solutions (3), we should find that contains the maximum number of positive $\lambda_j$. To do so, it suffices to solve the following MILP problem:

$$
\begin{aligned}
Max \quad & \sum_{j \in J} I_j \\
s.t. \quad & \sum_{j \in J} \lambda_j x_{ij} \leq \theta^* x_{ro}, \quad i = 1, \dots, m, \\
& \sum_{j \in J} \lambda_j y_{rj} \geq y_{ro}, \quad r = 1, \dots, s, \\
& \sum_{j \in J} \lambda_j = 1, \\
& I_j \leq M \lambda_j, I_j \in \{0,1\}, \lambda_j \geq 0, j \in J,
\end{aligned}
\quad (6)
$$

where $\theta^*$ is the optimal objective of (3); and $M$ is a sufficiently large positive quantity.

The idea behind the formulation of Model (6) is as follows: note that $I_j^* = 1$ if and only if $\lambda_j^* > 0$, i.e., $DMU_j$ enters actively into the optimal evaluation by (3). Since we are maximizing $\sum_{j \in J} I_j$ and $I_j \in \{0,1\}$, Model (6) will be directed towards finding optimal solutions with as many $I_j^* = 1$ as possible, i.e., with as many possible $\lambda_j^* > 0$; equivalently, with as many reference DMUs that can be entered into the optimal evaluation by (3).

Therefore, solving (6) yields such an optimal solution to (3) whose corresponding RS contains the maximum number of the reference DMUs to $\text{DMU}_o$. To be precise, we have the following theorem.

**Theorem 1.** Suppose that $\lambda^*$ is an optimal solution to Model (6). Then, $R_o^M = \{\text{DMU}_j | \lambda_j^* > 0\}$.

**Proof.** Let $\lambda^*$ is an optimal solution to Model (6) with the optimal objective $I_o^*$ and $S_o = \{\text{DMU}_j | \lambda_j^* > 0\}$. Then, $\lambda^*$ is an optimal solution to (3). By Definitions 2 and 3, we conclude that $S_o \subseteq R_o^M$. Now, if $R_o^M \subseteq S_o$, we are done; otherwise, there must exist a reference DMU such as $\text{DMU}_k$ which belongs to $R_o^M$ but not to $S_o$, i.e., $\lambda_k^* = 0$. Therefore, by Definition 3, there must exist an optimal solution to (3) such as $\hat{\lambda}$ where $\hat{\lambda}_k > 0$. By choosing a strict convex combination of the optimal solutions $\lambda^*$ and $\hat{\lambda}$, we can construct an optimal solution to (3) that is a feasible solution to (6) with the objective value greater than $I_o^*$; this is a contradiction. ∎

### 3.1.1. Computational Procedure for Solving Model (6)

Since Model (6) is a binary MILP problem, it is hard to solve directly using MILP solvers. So, it is imperative to enhance the computational capability of this model. To do so, integrating the LP and DEA concepts, we propose an LP-based procedure for solving Model (6): Let $(\theta^*, \lambda^*, S^{-*}, S^{+*})$ is a basic optimal solution to (3), which is obtained via the Simplex method. Then, we consider two cases:

➢ The relative cost coefficient of all the non-basic variables, denoted by $z_j - c_j$,[1] are less than zero. Then, the solution is the unique optimal solution to (3), and $R_o^M = \{\text{DMU}_j | \lambda_j^* > 0\}$, terminate[2]. This is precisely because of the *complementary slackness theorem* of LP - if the relative cost coefficient of a non-basic variable be less than zero in the optimality, then this variable's value must be equal to zero in every optimum solution.

➢ The relative cost coefficient of all non-basic variables are not less than zero. Then, we consider three sets $R_o = \{\text{DMU}_j | \lambda_j^* > 0\}$, $S_o = \{\text{DMU}_j | \lambda_j^* \geq 0, z_j - c_j < 0\}$ and $T_o = \{\text{DMU}_j | \lambda_j^* \geq 0, z_j - c_j = 0\}$. As already mentioned, the members of $S_o$ cannot enter actively into the optimal evaluation by (3), and are zero in every optimum solution. Thus, we omit them and their coefficients in (3), and formulate the following LP to identify such elements of $T_o$ that can be the possible candidates for reference DMUs:

---

[1] For more details and used notations see Murty [16].
[2] Specifically, if in an optimal solution all basic $\lambda_j^*$ variables are positive and the relative cost coefficient of all nonbasic $\lambda_j^*$ variables are less than zero (<0), then $R_o^M = \{\text{DMU}_j | \lambda_j^* > 0\}$, terminate.

$$Max \quad \gamma_o = \sum_{j \in T_o} \lambda_j$$
$$s.t. \quad \sum_{j \in J-S_o} \lambda_j x_{ij} \leq \theta^* x_{io}, \quad i = 1, \dots, m,$$
$$\sum_{j \in J-S_o} \lambda_j y_{rj} \geq y_{ro}, \quad r = 1, \dots, s, \qquad (7)$$
$$\sum_{j \in J-S_o} \lambda_j = 1,$$
$$\lambda_j \geq 0, j \in J - S_o.$$

Having solved (7), we distinguish two cases:

**Case (I).** $\gamma_o^* = 0$: In this case, the set $R_o$ does not change and we have $R_o = R_o^M$, terminate.

**Case (II).** $\gamma_o^* > 0$: Let $\bar{\lambda}$ is an optimal solution to (7) and let $L_o = \{DMU_j | \bar{\lambda}_j > 0\} \subseteq T_o$ and $K_o = \{DMU_j | \bar{\lambda}_j \geq 0, z_j - c_j < 0$ with respect to the objective function of $(3)\} \subseteq T_o$. We emphasize that the optimal base obtained via solving (7) is also an optimal base for (3) and the relative cost coefficient in $K_o$ is computed with respect to the cost coefficient of the objective (3). Then, we update the new sets $R_o \leftarrow R_o \cup L_o$, $S_o \leftarrow S_o \cup K_o$ and $T_o \leftarrow T_o - (L_o \cup K_o)$.

We repeat this routine and solve (7) again with the updated index sets, $R_o$, $S_o$ and $T_o$. As the cardinal number of the set $T_o$ decreases by one or more in each step, the procedure should terminate in at most q (q is the cardinal number of $T_o$) steps.

## 3.2. Dual Approach

The dual of Model (3), which is called the *multiplier* form of the BCC model, is as follows:

$$Max \quad U^t y_o + u_o$$
$$s.t. \quad V^t x_o = 1,$$
$$U^t y_j - V^t x_j + u_o + t_j = 0, j \in J, \qquad (8)$$
$$U \geq 0, V \geq 0, t_j \geq 0, j \in J,$$

where $(V, U, u_o)$ is the vector of dual variables, and $V \in \mathbb{R}_{\geq 0}^m$, $U \in \mathbb{R}_{\geq}^s$; $u_o$ is a scalar variable associated with the convexity constraint, whose sign is unrestricted.

Let $(U^*, V^*, u_o^*, t^*)$ is an optimal solution to Model (8), when $DMU_o$ is evaluated. Then, $H: U^{*t} y - V^{*t} x + u_o^* = 0$ would be a *supporting hyperplane* of the PPS, which all the members of $R_o^M$ lie on this hyperplane [12]. Note that Model (8) may have alternative optimal solutions; each defines a supporting hyperplane passing through $DMU_o$. From the optimal solutions to (8), we aim to select that associated with the hyperplane passing *only* through the members of $R_o^M$.

First, note that the variables $\lambda_j$ and $t_j$, $j \in J$, in Models (3) and (8) are the complementary slackness variables. Now recall that to determine $R_o^M$ via Model (3), an optimal solution with the maximum number of $\lambda_j^* > 0$ should be found. Thus, considering the complementary slackness conditions, we should find a dual optimal solution with the minimum number of $t_j^* = 0$; i.e., we should find a dual optimal solution with the minimum number of DMUs lying on its associated hyperplane. Based on this observation, we design the following MILP formulation to reach our aim:

$$\begin{aligned}
Max \ I_o &= \sum_{j \in J} I_j \\
s.t. \ U^t y_o + u_o &= \theta^* \\
V^t x_o &= 1, \\
U^t y_j - V^t x_j + u_o + t_j &= 0, \ j \in J, \\
I_j &\leq M t_j, \ j \in J, \\
U \geq 0, V \geq 0, I_j &\in \{0,1\}, t_j \geq 0, \ j \in J.
\end{aligned} \tag{9}$$

As a visual illustration, consider $DMU_7$ in Figure 1. Obviously, when $DMU_7$ is evaluated via Model (8), the MRS is the singleton set $R_7^M = \{DMU_2\}$. In this evaluation, (8) has alternative optimal solutions that define an infinite number of supporting hyperplanes passing through $DMU_7$, including $H$ and the hayperplanes passing through the line segments $DMU_1 - DMU_2$ and $DMU_2 - DMU_4$. It can be seen that $DMU_1$ and $DMU_4$, respectively lying on the two latter hyperplanes, do not belong to $R_7^M$. However, $DMU_2$ lying on the hyperplane $H$ is the only member of $R_7^M$, lies on the hyperplane $H$. Model (9) yields such a supporting hyperplane $H$, which is binding only at the members of $R_7^M$ lie on it.

**Theorem 2.** Let $(U^*, V^*, u_o^*, t^*)$ is an optimal solution to Model (9), and define $J_o = \{DMU_j | t_j^* = 0\}$. Then $J_o$ is equal to $R_o^M$.

**Proof.** It is clear that $R_o^M \subseteq J_o$. Without loss of generality, suppose that $J_o = \{1, \dots, k\}$; so the optimal value of $I_o$ is $I_o^* = |J| - k = n - k$. If $J_o \subseteq R_o^M$, then we are done. Otherwise, without loss of generality, we assume that $J_o - R_o^M = \{1, \dots, r\}$; Then, the values of the variables $\lambda_1, \dots, \lambda_r$ are zero in every optimal solution of Model (3). From the strong complementary slackness conditions, there exists some optimal solution $(\widehat{U}, \widehat{V}, \widehat{u}_o, \widehat{t})$ of Model (8), such that the values of the variables $t_1, \dots, t_r$ are positive. Thus, we have the following relations:

$$\begin{cases} U^{*t} y_j - V^{*t} x_j + u_o^* = 0, j \in J_o, (t_j^* = 0), \\ U^{*t} y_j - V^{*t} x_j + u_o^* < 0, j \in J - J_o, (t_j^* > 0), \end{cases} \tag{10}$$

and

$$\begin{cases} \hat{U}^t y_j - \hat{V}^t x_j + \hat{u}_o = 0, j \in R_o^M, (\hat{t}_j = 0), \\ \hat{U}^t y_j - \hat{V}^t x_j + \hat{u}_o < 0, j \in J_o - R_o^M, (\hat{t}_j > 0), \\ \hat{U}^t y_j - \hat{V}^t x_j + \hat{u}_o \leq 0, j \in J - J_o. \end{cases} \quad (11)$$

By choosing a strict convex combination of the optimal solutions $(U^*, V^*, u_o^*, t^*)$ and $(\hat{U}, \hat{V}, \hat{u}_o, \hat{t})$, we obtain an optimal solution $(\bar{U}, \bar{V}, \bar{u}_o, \bar{t})$ to Model (8) which satisfies the following relations:

$$\begin{cases} \bar{U}^t y_j - \bar{V}^t x_j + \bar{u}_o = 0, j \in R_o^M, (\bar{t}_j = 0), \\ \bar{U}^t y_j - \bar{V}^t x_j + \bar{u}_o < 0, j \in J_o - R_o^M, (\bar{t}_j > 0). \end{cases} \quad (12)$$

Obviously, $(\bar{U}, \bar{V}, \bar{u}_o, \bar{t})$ is a feasible solution of Model (9) with the objective value $\bar{I}_o = |J| - (k - r) = n - k + r > n - k = I_o^*$; this is a contradiction. ∎

**Remark 2.** It is worth noting that the algorithm proposed by Jahanshahloo et al. [12] to determine the MRS based on (8), can be used as an iterative computational treatment for the MILP (9).

### *3.3. Extension to Other DEA Models*

As we know, the DEA models are generally classified into two diverse groups: radial and non-radial models. Our approach can be readily used for the group of radial DEA models without any further consideration, since the MRS of any given DMU is unique - due to the uniqueness of the projection point. In spite of possibility of multiple projection points concerning the non-radial DEA models, our approach can be also implemented for this group. To facilitate the development, we divide the non-radial DEA models into two categories:

**Category (i).** The models with linear objective functions such as the additive [7] and RAM models [8]:

In each model of this category, as demonstrated in Krivonozhko et al. [14], there exists a minimal face containing all projection points to the evaluated DMU. To determine all reference DMUs corresponding to all projection points it suffices to fix the optimal objective and maximize the number of positive intensity variables $\lambda_j$ as a second phase. For instance, let $\sigma_o^*$ is the optimal objective of the additive model. Then, the DMUs with positive $\lambda_j^*$ in the following MILP are the all reference DMUs to a given DMU$_o$:

$$\begin{aligned}
Max \quad & \sum_{j \in J} I_j \\
s.t. \quad & \sum_{i=1}^{m} s_i^- + \sum_{r=1}^{s} s_r^+ = \sigma_o^* \\
& \sum_{j \in J} \lambda_j x_{ij} + s_i^- = x_{io}, \quad i = 1, \ldots, m, \\
& \sum_{j \in J} \lambda_j y_{rj} - s_r^+ = y_{ro}, \quad r = 1, \ldots, s, \\
& \sum_{j \in J} \lambda_j = 1, \\
& I_j \leq M\lambda_j, I_j \in \{0,1\}, \lambda_j \geq 0, j \in J, \\
& s_i^- \geq 0, s_r^+ \geq 0, \forall i, r.
\end{aligned} \quad (13)$$

**Category (ii).** The models with non-linear objective functions such as the Russell and ERM models [17,18]:

In each model of this category, as one of the anonymous reviewers pointed out, there may exist multiple projection points which are not on the same facet. To identify all the reference DMUs, we need to identify the MRS pertaining to each projection point. Regarding the ERM model, we describe how to identify the MRS for each projection point.

Suppose that $(\lambda_j^*, \theta_i^*, \varphi_r^*, \forall j, i, r)$ is an optimal solution to the ERM model. To determine the MRS for the projection point $P = (\theta_1^* x_{1o}, \ldots, \theta_m^* x_{mo}, \varphi_1^* y_{1o}, \ldots, \varphi_r^* y_{so})$, we solve the following MILP:

$$\begin{aligned}
Max \quad & \sum_{j \in J} I_j \\
s.t. \quad & \sum_{j \in J} \lambda_j x_{ij} = \theta_{io}^* x_{io}, \quad i = 1, \ldots, m, \\
& \sum_{j \in J} \lambda_j y_{rj} = \varphi_{ro}^* y_{ro}, \quad r = 1, \ldots, s, \\
& \sum_{j \in J} \lambda_j = 1, \\
& I_j \leq M\lambda_j, I_j \in \{0,1\}, \lambda_j \geq 0, j \in J, \\
& \theta_{io} \geq 0, \varphi_{ro} \geq 0, \forall i, r.
\end{aligned} \quad (14)$$

Having identified all the MRS corresponding to all projection points, finally, we collect them into a general set containing all reference DMUs.

As an example of the existence of multiple projection points in non-radial DEA models, consider $DMU_7$ in Figure 1. In the evaluation of $DMU_7$ via the RAM model, there are multiple projection points, $DMU_2, DMU_3,$ and $DMU_4$, which are on a same facet and $R_7^M = \{DMU_2, DMU_3, DMU_4\}$. Moreover, all of the reference DMUs are Pareto-efficient i.e., $E_7^M = \{DMU_2, DMU_3, DMU_4\}$[3].

It is worth noting that, Similar to the treatment prescribed for (6), a computational procedure can be adapted for (13) and (15).

## 4. Conclusion

It is well known that the occurrence of multiple reference sets is a crucial issue is the identification of all reference DMUs to a given DMU. To resolve this difficulty, we first introduced the useful notion of the MRS containing all the reference DMUs. Then with regard to both primal and dual formulations, we proposed a MILP-based approach to determine the MRS. Since this approach is compatible with both primal and dual forms of the DEA models, it can be very useful for real-world applications in which the number of DMUs is large. In such applications, the use of primal-based MILP would be easier and computationally efficient than its corresponding dual formulation because of having less number of constraints. However, there may exist some cases in which the use of dual-based MILP may be advantageous because of technical purposes. Our approach is more general than the existing ones as it is applicable to all kind of - radial and non-radial - DEA models.

In the proposed approach, the problem of finding MRS is formulated as a unique MILP. To treat effectively with the primal-based MILP, we established an LP-based computational procedure, and demonstrated that the algorithm by Jahanshahloo et al. [12] can be used as a computational treatment for our dual-based MILP.

We used a simple numerical example to describe our approach and its properties throughout the paper; however, one can use this approach for any real-world application.

---

[3] Due to taking all non-zero slacks into account in the non-radial DEA, it can be easily proved that all reference DMUs are Pareto-efficient (see [9]) and thus $R_o^M = E_o^M$.